**Title** Solving polynomial differential equations by transforming them to linear functional-differential equations

John Michael Nahay
Broadway Performance Systems, LLC
1013 Rosemere Ave., Silver Spring, MD 20904-3008

Email: resolvent@comcast.net

**Abstract**

We present a new approach to solving polynomial ordinary differential equations by transforming them to linear functional equations and then solving the linear functional equations. We will focus most of our attention upon the first-order Abel differential equation with two nonlinear terms in order to demonstrate in as much detail as possible the computations necessary for a complete solution. We mention in our section on further developments that the basic transformation idea can be generalized to apply to differential equations of any order, to a system of ordinary differential equations without first differentially eliminating the multiple dependent variables, and even to partial differential equations.

**Notation**

For each positive integer $K$ define $[K]$ to be the set of integers $\{k \ni 1 \le k \le K\}$.

**Introduction**

This approach is different than Yu N Kosovtsov's chronological operator algebra method [1]. We will first focus our attention upon the nonlinear first-order Abel [9] differential equation

$$\left(\frac{dz}{dx}\right)^n = \sum_{k=1}^{K} g_k(x) \cdot z^{m_k} \qquad (1)$$

where $x$ is the independent variable and $z$ is the dependent variable. We use the letter $g_k(x)$ to denote a *given* sufficiently differentiable function of $x$. We make the key hypothesis that the dependence of $z$ upon $x$ occurs only through the $g_k(x)$ and that this dependence is *continuously differentiable.* Without loss of generality, we assume that any solution $z$ of (1) can be expressed as a *multirariable function*

$$z = F(g_1(x), ..., g_K(x), m_1, ..., m_K, n) \qquad (2)$$

of $2K+1$ variables – or slots. We suppress the dependence of $z$ upon the arbitrary constant of integration. We assume that this multivariable function is analytic at 0 around the first set of $K$ slots and analytic at $n$ around the second set of $K$ slots. In other words, we assume there exists some multivariable function

$$F(u_1,...,u_K,m_1,...,m_K,n) \qquad (3)$$

of $2K+1$ indeterminate variables $\{u_1,...,u_K,m_1,...,m_K,n\}$ such that when the indeterminate $u_k$ is replaced by $g_k(x)$ for each $k\in[K]$ we obtain a solution of (1). There is no confusion in identifying $m_k$ and $n$ as both an indeterminate in function (3) and as the particular complex- or real-valued number given in (1).

**The transformation**

The multivariable function in (3) is *not* unique in general, because we replace $K$ algebraically independent indeterminates $u_k$ with $K$ functions $g_k(x)$ which are all related to the same variable $x$. We will discuss this more on the section titled "Non-uniquess of the solution". Let $\alpha$ be an indeterminate constant with respect to $x$. Multiply equation (1) by $z^\alpha$. We obtain

$$z^\alpha\cdot\left(\frac{dz}{dx}\right)^n=\sum_{k=1}^{K}g_k(x)\cdot z^{m_k+\alpha} \qquad (4)$$

Rewrite (4) as

$$\left(z^{\alpha/n}\cdot\frac{dz}{dx}\right)^n=\sum_{k=1}^{K}g_k(x)\cdot z^{m_k+\alpha} \qquad (5)$$

Define

$$w\equiv z^{\alpha/n+1}=z^{\frac{n+\alpha}{n}}. \qquad (6)$$

So

$$z=w^{\frac{n}{n+\alpha}} \qquad (7)$$

and

$$\frac{dw}{dx}=\frac{n+\alpha}{n}\cdot z^{\frac{\alpha}{n}}\cdot\frac{dz}{dx}. \qquad (8)$$

So equation (5) becomes

$$\left(\frac{dw}{dx}\right)^n=\sum_{k=1}^{K}\left(\frac{n+\alpha}{n}\right)^n\cdot g_k(x)\cdot w^{\frac{n\cdot(m_k+\alpha)}{n+\alpha}}. \qquad (9)$$

Equation (9) shows that $w$ satisfies the same differential equation as (1) but with $g_k(x)$ replaced with $\left(\frac{n+\alpha}{n}\right)^n\cdot g_k(x)$ and $m_k$ replaced with $n\cdot\frac{m_k+\alpha}{n+\alpha}$. Hence, $u$ must be given by the same functional form as (3). Hence

$$w=F(\left(\frac{n+\alpha}{n}\right)^n\cdot g_1(x),...,\left(\frac{n+\alpha}{n}\right)^n\cdot g_K(x),\frac{n\cdot(m_1+\alpha)}{n+\alpha},...,\frac{n\cdot(m_K+\alpha)}{n+\alpha},n) \qquad (10)$$

But (6) relates $F$ before this functional substitution to $F$ after this functional substitution. In other words, $F$ satisfies the functional, *non-differential* relation

$$\left(F\left(\left(\frac{n+\alpha}{n}\right)^n\cdot g_1(x),...,\left(\frac{n+\alpha}{n}\right)^n\cdot g_K(x),\frac{n\cdot(m_1+\alpha)}{n+\alpha},...,\frac{n\cdot(m_K+\alpha)}{n+\alpha},n\right)\right)^{\frac{n}{n+\alpha}} \tag{11}$$
$$= F(g_1(x),...,g_K(x),m_1,...,m_K,n)$$

We have *not* yet associated $m_k$ with $g_k(x)$. Observe that one may permute the functions and corresponding powers in (11). This association will be made when we fix *initial* conditions on the $m_k$ and $g_k(x)$, that is, when we solve for $F(u_1,...,u_K,m_1,...,m_K,n)$ in terms of $F(u_1,...,u_{K-1},0,m_1,...,m_{K-1},1,n)$, and so forth.

Note also that if $z = F(g_1(x),...,g_K(x),m_1,...,m_K,n)$ satisfies (11), then any power of it, $H(g_1(x),...,g_K(x),m_1,...,m_K,n) \equiv (F(g_1(x),...,g_K(x),m_1,...,m_K,n))^\beta$, satisfies the same nonlinear functional equation, because

$$\left(H\left(\left(\frac{n+\alpha}{n}\right)^n\cdot g_1(x),...,\left(\frac{n+\alpha}{n}\right)^n\cdot g_K(x),\frac{n\cdot(m_1+\alpha)}{n+\alpha},...,\frac{n\cdot(m_K+\alpha)}{n+\alpha},n\right)\right)^{\frac{n}{n+\alpha}}$$

$$\left(F\left(\left(\frac{n+\alpha}{n}\right)^n\cdot g_1(x),...,\left(\frac{n+\alpha}{n}\right)^n\cdot g_K(x),\frac{n\cdot(m_1+\alpha)}{n+\alpha},...,\frac{n\cdot(m_K+\alpha)}{n+\alpha},n\right)\right)^{\beta\cdot\frac{n}{n+\alpha}} \tag{12}$$
$$= (F(g_1(x),...,g_K(x),m_1,...,m_K,n))^\beta$$
$$= H(g_1(x),...,g_K(x),m_1,...,m_K,n)$$

Define

$$\bar{u}_k \equiv \left(\frac{n+\alpha}{n}\right)^n\cdot u_k.\ \bar{g}_k(x) \equiv \left(\frac{n+\alpha}{n}\right)^n\cdot g_k(x), \text{ and } \bar{m}_k \equiv \frac{n\cdot(m_k+\alpha)}{n+\alpha}. \tag{13}$$

When we wish to emphasize a substitution or "transformation" from $u_k$ to $\bar{u}_k$, we will sometimes write $u_k \to \bar{u}_k$. Similarly, when we subsitute $g_k(x)$ for $u_k$, we will write $u_k \to g_k(x)$. When we substitute $\bar{m}_k$ for $m_k$, we will write $m_k \to \bar{m}_k$.

In order to shorter our notation even more, define

$$F(u,m) \equiv F(u_1,...,u_K,m_1,...,m_K) \tag{14}$$

and $$F(\bar{u},\bar{m}) \equiv F(\bar{u}_1,...,\bar{u}_K,\bar{m}_1,...,\bar{m}_K) \tag{15}$$

and $$F(g,m) \equiv F(g_1(x),...,g_K(x),m_1,...,m_K) \tag{16}$$

and $$F(\bar{g},\bar{m}) \equiv F(\bar{g}_1(x),...,\bar{g}_K(x),\bar{m}_1,...,\bar{m}_K) \tag{17}$$

and $$\bar{F} \equiv F(\bar{u},\bar{m}) \tag{18}$$

So (11) becomes

$$F(g,m) = (F(\bar{g},\bar{m}))^{\frac{n}{n+\alpha}} \tag{19}$$

We dropped showing the dependence upon $n$ because $n$ does not get transformed in (11). Take the natural logarithm of (19)

$$\ln(F(g,m)) = \frac{n}{n+\alpha}\ln(F(\bar{g},\bar{m})) \tag{20}$$

Observe that if there exists a function (3), which satisfies the same functional form as (20), in other words, if there exists a function, which satisfies

$$\ln(F(u,m)) = \frac{n}{n+\alpha}\ln(F(\overline{u},\overline{m})) \qquad (21)$$

then the function will satisfy (20) when $g_k(x)$ is substituted for $u_k$.
Define

$$F_{u,k} \qquad (22)$$

to be the partial derivative of $F(u_1,...,u_K,m_1,...,m_K)$ with respect to $u_k$, holding all the $u_{j\neq k}$ and all the $m_j$ fixed.
Define

$$F_{m,k} \qquad (22)$$

to be the partial derivative of $F(u_1,...,u_K,m_1,...,m_K)$ with respect to $m_k$, holding all the $m_{j\neq k}$ and all the $u_j$ fixed.
Define

$$\overline{F}_{\overline{u},k} \qquad (22)$$

to be the partial derivative of $F(\overline{u}_1,...,\overline{u}_K,\overline{m}_1,...,\overline{m}_K)$ with respect to $\overline{u}_k$, holding all the $\overline{u}_{j\neq k}$ and all the $\overline{m}_j$ fixed.
Define

$$\overline{F}_{\overline{m},k} \qquad (22)$$

to be the partial derivative of $F(\overline{u}_1,...,\overline{u}_K,\overline{m}_1,...,\overline{m}_K)$ with respect to $\overline{m}_k$, holding all the $\overline{m}_{j\neq k}$ and all the $\overline{u}_j$ fixed.

Differentiate (21) with respect to each of the $u_k$. We obtain

$\frac{F_{u,k}}{F} = \frac{n}{n+\alpha}\frac{F_{\overline{u},k}}{\overline{F}}\cdot\frac{d\overline{u}_k}{du_k}$. From (13) we have $\frac{d\overline{u}_k}{du_k} = \left(\frac{n+\alpha}{n}\right)^n$. So

$$\frac{F_{u,k}}{F} = \left(\frac{n+\alpha}{n}\right)^{n-1}\frac{F_{\overline{u},k}}{\overline{F}} \qquad (23)$$

Define

$$\Lambda_{u,k} \equiv \frac{F_{u,k}}{F}. \qquad (24)$$

Then (23) states that, when the transformations $u_j \to \overline{u}_j$ and $m_j \to \overline{m}_j$ are all made, then the transformation $\Lambda_{u,k} \to \left(\frac{n+\alpha}{n}\right)^{n-1}\Lambda_{\overline{u},k}$ is made. In other words, $\Lambda_{u,k}$ satisfies the *linear functional equation*

$$\Lambda_{u,k} = \left(\frac{n+\alpha}{n}\right)^{n-1}\Lambda_{\overline{u},k} \qquad (25)$$

Differentiate (21) with respect to each of the $m_k$. We obtain $\frac{F_{m,k}}{F} = \frac{n}{n+\alpha}\frac{F_{\bar{m},k}}{\bar{F}} \cdot \frac{d\bar{m}_k}{dm_k}$. From (13) we have $\frac{d\bar{m}_k}{dm_k} \equiv \frac{n}{n+\alpha}$. So

$$\frac{F_{m,k}}{F} = \left(\frac{n}{n+\alpha}\right)^2 \frac{F_{\bar{m},k}}{\bar{F}} \tag{26}$$

Define

$$\Omega_{m,k} \equiv \frac{F_{m,k}}{F}. \tag{27}$$

Then (26) states that when the transformations $u_j \to \bar{u}_j$ and $m_j \to \bar{m}_j$ are all made then the transformation $\Omega_{m,k} \to \left(\frac{n}{n+\alpha}\right)^2 \Omega_{\bar{m},k}$ is made. In other words, $\Omega_{m,k}$ satisfies the *linear functional equation*

$$\Omega_{m,k} = \left(\frac{n}{n+\alpha}\right)^2 \Omega_{\bar{m},k} \tag{28}$$

So, we have reduced the problem of solving (1) to solving the linear functional equations (25) and (28). Observe that the transformation $m_j \to \bar{m}_j$ results in $m - n \to \frac{n(m+\alpha)}{n+\alpha} - n = \left(\frac{n}{n+\alpha}\right) \cdot (m-n)$. From this we see that the general solution of (25) is

$$\Lambda_{u,k} = \sum_I {}_{u,k}c_I \cdot \prod_{j=1}^{K} u_j^{a_{I,j}} \cdot \prod_{j=1}^{K} (m_j - n)^{b_{I,j}} \tag{29}$$

where the $a_{I,j}$ and the $b_{I,j}$ are subject to

$$n \cdot \sum_{j=1}^{K} a_{I,j} - \sum_{j=1}^{K} b_{I,j} = n - 1 \tag{30}$$

where the ${}_{u,k}c_I$ must be *constant* with respect to *all* $u_j$ and *all* $m_j$. We see that the general solution of (26) is

$$\Omega_{m,k} = \sum_I {}_{m,k}c_I \cdot \prod_{j=1}^{K} u_j^{a_{I,j}} \cdot \prod_{j=1}^{K} (m_j - n)^{b_{I,j}} \tag{31}$$

where the $a_{I,j}$ and the $b_{I,j}$ are subject to

$$n \cdot \sum_{j=1}^{K} a_{I,j} - \sum_{j=1}^{K} b_{I,j} = 2 \tag{32}$$

where the ${}_{u,k}c_I$ must be *constant* with respect to *all* $u_j$ and *all* $m_j$.

We assume that $F$ is analytic at 0 in each of the $u_k$ and analytic at $n$ in each of the $m_k$ separately. Hence, for each $k \in [K]$, $\Lambda_{u,k}$ and $\Omega_{m,k}$ are analytic at 0 in each of the $u_k$ and analytic at $n$ in each of the $m_k$ separately. This implies that all the $a_{I,j}$ and $b_{I,j}$ in (29) and (31) are nonnegative.

## Ladder of boundary conditions

When $n \neq 1$, it is not known if there is sufficient information – i.e. a sufficient number of relations - from (27) and (29) to fully solve (1). Nevertheless, we *do* have the ladder of boundary conditions which makes this current method of solution hopeful. Specifically, assume that, for each $j \in [K]$, $F(u_1,...,u_K,m_1,...,m_K)$ is known when $u_j = 0$ (or $u_j = u_{j'}$ for some $j' \neq j$) and that, for each $j \in [K]$, $F(u_1,...,u_K,m_1,...,m_K)$ is known when $m_j = n$ (or $m_j = 1$ or $m_j = m_{j'}$ for some $j' \neq j$). Then, these *known* functions constitute the *boundary* conditions for the linear functional equations (24) and (26). We can choose $m_j = m_{j'}$ for some $j' \neq j$ as a boundary condition because then those two terms of Abel's differential equation (1) can be combined with the same exponent.

The ladder of boundary conditions is created when we express $F(u_1,...,u_K,m_1,...,m_K)$, with some subset $S$ of the $u_k$ and/or $m_k$ specialized to the particular values suggested in the previous paragraph, in terms of $F(u_1,...,u_K,m_1,...,m_K)$ with some larger subset $S'$, with $S \subset S'$, of the $u_k$ and/or $m_k$ specialized to the particular values suggested in the previous paragraph.

When $n = 1$, we have an auxiliary linear mixed partial functional-differential equation which we can use to obtain a solution of (1).

## The linear mixed partial functional-differential equation for n=1

Define $F_{g,k}$ to be $F_{u,k}$ from definition (22) with each indeterminate $u_k$ replaced with the function $g_k(x)$. Elementary differentiation implies

$$\frac{dz}{dx} = \frac{dF}{dx} = \sum_{k=1}^{K} F_{g,k} \cdot \frac{dg_k(x)}{dx}. \tag{33}$$

When $n = 1$, (11) simplies to

$$\begin{aligned} &F((1+\alpha)\cdot g_1(x),...,(1+\alpha)\cdot g_K(x), \frac{m_1+\alpha}{1+\alpha},...,\frac{m_K+\alpha}{1+\alpha}) \\ &= (F(g_1(x),...,g_K(x),m_1,...,m_K))^{1+\alpha} \end{aligned} \tag{34}$$

Since $\alpha$ is arbitrary, we may substitute $m_k - 1$ for $\alpha$ in (34) to obtain

$$z^{m_k} = F^{m_k} = F(m_k \cdot g_1(x),...,m_k \cdot g_K(x), \frac{m_1+m_k-1}{m_k},...,\frac{m_K+m_k-1}{m_k}) \tag{35}$$

So $$\sum_{k=1}^{K} g_k(x)\cdot z^{m_k} = \sum_{k=1}^{K} g_k(x)\cdot F(m_k \cdot g_1(x),...,m_k \cdot g_K(x), \frac{m_1+m_k-1}{m_k},...,\frac{m_K+m_k-1}{m_k}) \tag{36}$$

The Abel differential equation is formed by equating (33) and (36). Hence, if we can find a function $F$ written as in (3) which satisfies

$$\sum_{k=1}^{K} F_{u,k} \cdot \frac{dg_k(x)}{dx} = \sum_{k=1}^{K} u_k \cdot F(m_k \cdot u_1,...,m_k \cdot u_K, \frac{m_1+m_k-1}{m_k},...,\frac{m_K+m_k-1}{m_k}) \tag{37}$$

in indeterminates $u_k$, then we will have found a function $F$ which satisfies Abel's differential equation when $u_k \to g_k(x)$.

**Non-uniqueness of the solution**

As written now, (37) is not quite correct. We would need to introduce a lot of notation which we will not need elsewhere, so it was deemed best to leave (37) as is. Specifically, the term $\frac{dg_k(x)}{dx}$ should not be in (37). We need to express $\frac{dg_k(x)}{dx}$ in terms of the *indeterminates* $u_k$. There is *no unique way to do this*. Suppose that each of the $g_k(x)$ is *theoretically invertible* on some restricted domain. Although not rigorously correct, since the $u_k$ are indeterminates and, hence, have no "attachment" to $x$, we will temporarily write for this section the inverse as $x = g_k^{(-1)}(u_k)$. But, since we assume this true for each $k \in [K]$, it follows that $x$, and hence any function $r(x)$ of $x$, can be expressed in infinitely ways as a multivariable function $F(u_1,...,u_K)$ of the $u_k$'s such that when the substitution $u_k \to g_k(x)$ is made in $F$ we recover $r(x) = F(g_1(x),...,g_K(x))$.

**Example.** Suppose $g_1(x) = \frac{1}{2+x}$ and $g_2(x) = x^2$. Suppose $r(x) \equiv \frac{dg_1(x)}{dx} = -(2+x)^{-2}$.

We have $x = g_1^{(-1)}(u_1) = \frac{1}{u_1} - 2$ and $x = g_2^{(-1)}(u_2) = \sqrt{u_2}$. So, we can write $r(x)$ in the form

$$r(x) = -\frac{1}{2+x}\cdot\frac{1}{2+x} = -u_1 \cdot \frac{1}{2+\sqrt{u_2}} \quad \text{so } F(u_1,u_2) = -u_1 \cdot \frac{1}{2+\sqrt{u_2}}$$

or

$$r(x) = -\frac{1}{2+x}\cdot\frac{1}{2+x} = -u_1 \cdot u_1 \quad \text{so } F(u_1,u_2) = -u_1^2$$

or

$$r(x) = -\frac{1}{2+x}\cdot\frac{1}{2+x} = -\frac{1}{2+\sqrt{u_2}}\cdot\frac{1}{2+\sqrt{u_2}} \quad \text{so } F(u_1,u_2) = -\left(\frac{1}{2+\sqrt{u_2}}\right)^2$$

The dilemma in solving the Abel differential equation by our current method is that if we are given a function such as $r(x) = g_1(x) + g_2(x)$, expressed in terms of the $g_k(x)$, we seek a "canonical" form, namely, $F(u_1,u_2) = u_1 + u_2$. This "canonical" solution will be the "most" symmetrical function of the $u_k$'s.

We will explore the technical difficulties of this non-uniqueness of a multivariable function of indeterminates on a simpler case of equation (1) – the two-term Abel differential equation.

**The two-term Abel differential equation**

Set $n=1$, $K=2$ in (1). Define $g(x)\equiv g_1(x)$, $h(x)\equiv g_2(x)$, $u\equiv u_1$, $v\equiv u_2$, $m\equiv m_1$, (redefine) $n\equiv m_2$ in (1) to get

$$\frac{dz}{dx}=g(x)\cdot z^m+h(x)\cdot z^n. \tag{38}$$

Define $\phi(x)\equiv\frac{dg(x)}{dx}$ and $\psi(x)\equiv\frac{dh(x)}{dx}$. At this time, we do not know in terms of which variables, $u$ and/or $v$, we wish to express $x$. So, temporarily write $x=f(u,v)$, knowing that in each place that the symbol $f(u,v)$ appears, it could be a *different* function of $u$ and/or $v$. Then (36) becomes

$$\begin{aligned}&F_u(u,v,m,n)\cdot\phi(f(u,v))+F_v(u,v,m,n)\cdot\psi(f(u,v))\\&=u\cdot F(m\cdot u,m\cdot v,\frac{m-1}{m}+1,\frac{n-1}{m}+1)+v\cdot F(n\cdot u,n\cdot v,\frac{m-1}{n}+1,\frac{n-1}{n}+1)\end{aligned} \tag{39}$$

Observe that Abel's equation (1) remains unchanged if we switch $g_j(x)$ with $g_k(x)$ as long as we switch $m_j$ with $m_k$ simultaneously. This symmetry in Abel's equation suggests a similar symmetry in the solution $F$, which suggests a similar symmetry in the linear functional-differential equation (39), which $F$ satisfies. Ideally, we want $\phi(f(u,v))$ to be a function of $u$ only, because $\phi(x)\equiv\frac{dg(x)}{dx}$ and $u$ is the variable that is replaced with $g(x)$. So, we choose $f(u,v)$ to be $g^{(-1)}(u)$ as the argument of $\phi$. By the same reasoning, we choose $f(u,v)$ to be $h^{(-1)}(v)$ as the argument of $\psi$. So (39) becomes

$$\begin{aligned}&F_u(u,v,m,n)\cdot\phi\circ g^{(-1)}(u)+F_v(u,v,m,n)\cdot\psi\circ h^{(-1)}(v)\\&=u\cdot F(m\cdot u,m\cdot v,\frac{m-1}{m}+1,\frac{n-1}{m}+1)+v\cdot F(n\cdot u,n\cdot v,\frac{m-1}{n}+1,\frac{n-1}{n}+1)\end{aligned} \tag{40}$$

Definitions (24) and (27) and the equality of mixed second-order partial derivatives $\Lambda_{uv}=\Lambda_{vu}$, $\Lambda_{um}=\Omega_{mu}$, $\Lambda_{un}=\Omega_{nu}$, $\Lambda_{vm}=\Omega_{mv}$, $\Lambda_{vn}=\Omega_{nv}$, $\Omega_{mn}=\Omega_{nm}$ imply the following relations

$$\begin{aligned}F(u,v,m,n)&=F(0,v,m,n)\cdot\exp\left(\int_0^u\Lambda_u(\theta,v,m,n)\cdot d\theta\right)\\&=F(u,0,m,n)\cdot\exp\left(\int_0^v\Lambda_v(u,\theta,m,n)\cdot d\theta\right)\\&=F(u,v,1,n)\cdot\exp\left(\int_0^m\Omega_m(u,v,\theta,n)\cdot d\theta\right)\\&=F(u,v,m,1)\cdot\exp\left(\int_0^n\Omega_n(u,v,m,\theta)\cdot d\theta\right)\end{aligned} \tag{41}$$

So, we have reduced the solution of the Abel differential equation (1) with $n = 1$ to the solution of a linear mixed partial functional-differential equation (37) with the assistance of auxiliary equations like those shown in (41) for the two-term Abel differential equation (38). The linear mixed partial functional-differential equation for the two-term Abel differential equation (38) is given by (40). We must use the ladder of boundary conditions suggested earlier to get the complete solution. For the two-term Abel differential equation (38), this ladder of boundary conditions is a ladder with two rungs.

**The Bernoulli Equation**

The Bernoulli ordinary differential equation is a special case of the two-term Abel ordinary differential equation. Specialize $n = 1$ in (38). Then

$$\frac{dz}{dx} = g(x)\cdot z^m + h(x)\cdot z\,. \tag{42}$$

The famous solution of (42) is

$$z = \left[ z_0 + (1-m)\int_{t=0}^{t=x} h(t)\cdot \exp\left( (m-1)\int_{\theta=0}^{\theta=t} g(\theta)d\theta \right) dt \right]^{1/(1-m)} \cdot \exp\left( \int_{\theta=0}^{\theta=x} g(\theta)d\theta \right) \tag{43}$$

In (41) we must somehow associate $u$ with $m$ and $v$ with $n$. We know what we want the properties of the functions $F(u,v,m,1)$ and $F(u,v,1,n)$ to be. We want $F(u,v,m,1)$ to be such that, upon specialization $u \to g(x)$ and $v \to h(x)$, $F(g(x),h(x),m,1)$ is (43). Similarly, $F(g(x),h(x),1,n)$ must reduce to (43) with $n$ in place of $m$ and the roles of $g(x)$ and $h(x)$ switched.

Similarly, $F(0,v,m,n)$ and $F(u,0,m,n)$ "collapse" to solutions of easy cases of the Abel equation. We know $F(0,v,m,n)$ is the solution of

$$\frac{dz}{dx} = v\cdot z^n\,. \tag{44}$$

So

$$F(0,v,m,n) = \left( (1-n)\cdot \int_{t=0}^{t=x} v(t)dt + z_0^{1-n} \right)^{1/(1-n)} \quad \text{for } n \neq 1 \tag{45}$$

and

$$F(0,v,m,1) = z_0 \cdot \exp\left( \int_{t=0}^{t=x} v(t)dt \right). \tag{46}$$

**Further developments**

Research has proceeded in two separate directions at the time of this writing. First, and most importantly, attempts are being made at computing the solution of the linear mixed partial functional-differential equation (40) for the two-term Abel differential equation (38). One of the great difficulties is expressing the boundary condition (43), when the two-term Abel equation reduces further to the Bernoulli equation (42), as a function $F(u,v,m,1)$ of indeterminates $u$ and $v$ in a canonical way, in a way that makes the solution of (40) as easy as possible. Computation of the solution $F(u,v,m,n)$ of (40) as a power series in the four variables $u$, $v$, $m$, and $n$ has been attempted on the computer algebra system Maple. The latest computations determined about $5^4 = 625$ terms, up to 4th degree in each of these 4 variables, of the solution in terms of the Bernoulli functions (43) and $F(u,v,1,n)$ and the solution (45) of the one-term Abel equation (44) and $F(0,v,m,n)$, when (43) and (45) were expressed as power series in $u$ and/or $v$ and/or $m$ and/or $n$. However, no general pattern could be ascertained. The most speculative idea is to relate linear functional-differential equations like (40) and equations like (41) to the author's earlier work on differential resolvents [3], [4], [5], [6], [7]. Differential resolvents have proven useful before. The author's work is cited in [1] and [8].

The second direction of current research is to extend the basic transformation idea (4) of multiplying the original differential equation by $z^{\alpha}$ to systems of polynomial partial differential equations. Initial attempts look extremely hopeful. The generalization to systems of partial differential equations proceeds by multiplying in products of arbitrary powers of *all* the partial derivatives of *all* orders of *all* the dependent variables which appear in the given equations. One obtains generalizations of (11), (25) and (28), which are *much* more messy and complicated, and bear *many* more terms, than their first-order scalar Abel ordinary differential equation counterpart. Furthermore, the corresponding counterpart to (37) for vector polynomial partial differential equations has not been discovered. It is not known whether such a corresponding mixed partial functional-differential equation is absolutely necessary, but it is suspected to be.

**Acknowledgements**

The author thanks Dr Paul Nahay, President of Broadway Performance Systems, LLC for his support of this research project as well as for his assistance in Java programming. The participants on the Mapleprimes webboard were also of great assistance in helping the author write code.

**References**

[1] Bostan, Alin, Frédéric Chyzak, Grégoire Lecerf, Bruno Salvy, & Éric Schost "Differential equations for algebraic functions" Proceedings of the International Symposium of Symbolic & Algebraic Computation (ISSAC 2007) (2007 March 23)

[2] Kosovtsov Yu N "The chronological operator algebra and formal solutions of differential equations" *arXiv:math-ph/0409035 v1 16 Sep 2004/* (2004 September 16)

[3] Nahay, John Michael "Joint differential resolvent for pseudopolynomials" Arxiv.org ID 0803.2477 and www.sci.ccny.cuny.edu/~ksda/PostedPapers/jn-mar-01-08.pdf (2008 March 1)

[4] Nahay, John Michael "Linear relations among algebraic solutions of differential equations" Journal of Differential Equations, Volume 191 Issue 2 pp. 323-347

(2003 July 1)

[5] Nahay, John Michael "Powersum formula for polynomials whose distinct roots are differentially independent over constants" International Journal of Mathematics and Mathematical Sciences, Volume 32 Issue 12 pp. 721-738 (2002 December 22)

[6] Nahay, John Michael "Differential resolvents of minimal order and weight" International Journal of Mathematics and Mathematical Sciences, Volume 2004 Issue 54 pp. 2867-2894 (2004 September 26)

[7] Nahay, John Michael "A partial factorization of the powersum formula" International Journal of Mathematics and Mathematical Sciences, Volume 2004 Issue 58 pp. 3075-3102 (2004 October 16)

[8] Visagie, Stephan E 2008 "On the efficient solvability of a simple class of nonlinear knapsack problems", Orion, Volume 24, No. 1, 2008 pp1-15.

[9] Zwillinger, D. *Handbook of Differential Equations, 3rd ed.* Boston, MA: Academic Press, p120 (1997)